\documentstyle[12pt]{article}
\newcommand{\bee}{\begin{eqnarray}}
\newcommand{\eee}{\end{eqnarray}}
\newcommand{\be}{\begin{eqnarray}}
\newcommand{\ee}{\end{eqnarray}}
\newcommand{\di}{\frac {\partial} {\partial x_i}}
\newcommand\ie{{\it i.e.}}
\newcommand\nn{\nonumber \\}
\newcommand\defeq{\stackrel {def}{=}}
\newcommand\half{\frac 1 2}

\newcommand{\R}{{\bf R}}
\newcommand{\x}{{\vec x}}
\newcommand{\vv}{{\vec v}}

\font\frtnfr=eufm10   scaled\magstep1
\font\twlfr=eufm10
\font\tenfr=eufm10

\newfam\frfam
\textfont\frfam=\frtnfr
\scriptfont\frfam=\twlfr
\scriptscriptfont\frfam=\tenfr
\def\fr{\fam\frfam}

\font\frtnopen=msbm10  scaled\magstep2
\font\twlopen=msbm10
\font\tenopen=msbm10

\newfam\openfam
\textfont\openfam=\frtnopen
\scriptfont\openfam=\twlopen
\scriptscriptfont\openfam=\tenopen
\def\open{\fam\openfam}

\font\frtnsf = cmss12 scaled\magstep1
\font\twlsf = cmss10
\font\tensf = cmss9

\newfam\Scfam
\textfont\Scfam = \frtnsf
\scriptfont\Scfam = \twlsf
\scriptscriptfont\Scfam = \tensf

\begin{document}

%%%%%%%%%%%%%%%%%%%%%%%%%%%%%%%%%%%%%%%%%%%%%%%%%%%%%%%%%%%%%%%%%%%%%%%%%%%
%%%%%%%%%%%%%%      TITLE PAGE     %%%%%%%%%%%%%%%%%%%%%%%%%%%%%%%%%%%%%%%%
%%%%%%%%%%%%%%%%%%%%%%%%%%%%%%%%%%%%%%%%%%%%%%%%%%%%%%%%%%%%%%%%%%%%%%%%%%%

\sloppy
\title
 {
      \hfill{\normalsize\sf FIAN/TD/97-18}    \\
            \vspace{1cm}
      Rational Calogero models based on rank-2 root systems:
      supertraces on the superalgebras of observables
 }
\author
 {
   S.E.Konstein
          \thanks
             {E-mail: konstein@td.lpi.ac.ru}
          \thanks
             {This work is supported by the Russian Fund for Basic Research,
              Grants 96-02-17314 and 96-15-96463.}
  \\
               {\small \phantom{uuu}}
  \\
               {\it {\small} I.E.Tamm Department of Theoretical Physics,}
  \\
               {\it {\small} P. N. Lebedev Physical Institute,}
  \\
               {\it {\small} 117924, Leninsky Prospect 53, Moscow, Russia.}
 }
\date{ }
%--------------------------------------------------------------------------
\maketitle
%--------------------------------------------------------------------------
\begin{abstract}
\noindent
It is shown that the superalgebra $H_{W(I_2(n))}$ of observables of the
rational Calogero model based on the root system $I_2(n)$ possesses
$[(n+1)/2]$ supertraces. Model with three-particle
interaction based on the root system $G_2$ belongs
to this class of models and its superalgebra of observables has
3 independent supertraces.
\end{abstract}

\newpage
%%%%%%%%%%%%%%%%%%%%%%%%%%%%%%%%%%%%%%%%%%%%%%%%%%%%%%%%%%%%%%%%%%%%%%%%%%%
%%%%%%%%%%%%%%%%%%%%%%%%%%%%%%%%%%%%%%%%%%%%%%%%%%%%%%%%%%%%%%%%%%%%%%%%%%%

\bibliographystyle{nphys}
\renewcommand{\theequation}{\arabic{equation}}
\setcounter{equation}{0}
%---------------------------------------------

\section{Introduction}

In this work we continue to compute the number of supertraces
on the superalgebras of observables underlying the rational
Calogero models \cite{Cal} based on the root systems \cite{OP}.
The results for the root systems of $A_{N-1}$, $B_N$ ($C_N$) and $D_N$ types
are listed in \cite{KV}, \cite{K}. Here we consider the root systems
of $I_2(n)$ type.
%\cite{NB}.
As $W(A_2)=W(I_2(3))$,
$W(B_2)=W(I_2(4))$,
$W(D_2)=W(I_2(2))$ and $W(G_2)=W(I_2(6))$ the case under
consideration covers all rank-2 systems.
It is shown that the number of supertraces on the superalgebra
of observables of Calogero model based on this root system is $[(n+1)/2]$.
In particularly it gives that in the case of the root system $G_2$
there are 3 supertraces.

The definition and some properties of the superalgebra of observables
are discussed in the next section. In the section \ref{trace} the
condition sufficient for existence of the supertraces is formulated and some
consequences from the existence of the supertraces are discussed.
The superalgebra $H_{W(I_2(n))}$ of observables of the
rational Calogero model based on the root system $I_2(n)$
is described in the section \ref{I2} and the number of supertraces
on this superalgebra is computed in the last section.

\section{The superalgebra of observables}
The superalgebra of observables
of the rational Calogero model based on the root system
is defined in the following way.

Let us define the reflections
$R_\vv$ $\forall \vv \in {\open R}_N$, $\vv\neq 0$, as follows
\be\label{ref}
R_\vv (\x)=\x -2 \frac {(\x,\,\vv)} {(\vv,\,\vv)} \vv \qquad \forall \x
\in {\open R}_N.
\ee
Here the notation $(\cdot,\cdot)$ used for the inner product in ${\open R}_N$.
We will use also the coordinates of vectors, $v_i\defeq (\vv,\,\vec e_i)$,
where vectors $\vec e_i$ constitute the orthonormal basis in ${\open R}_N$:
$(\vec e_i,\, \vec e_j)=\delta_{ij}$.
The reflections (\ref{ref}) satisfy the properties
\be\label{prop}
R_\vv (\vv)=-\vv,\qquad R_\vv^2 =1,\qquad
({R}_\vv (\x),\,\vec u)=(\x,\,{R}_\vv (\vec u)),\quad
\forall \vv,\,\x,\,\vec u\in {\open R}_N.
\ee

The finite set of vectors $\R\subset {\open R}_N$ is the {\it root system}
if $\R$ is ${R}_\vv$-invariant $\forall \vv \in \R$ and
the group $W(\R)$ generated by all reflections ${R}_\vv$ with
$\vv \in \R$ (Coxeter group) is finite. The brief description of
the classification of the root systems one can find
for instance in \cite{OP}.

The group $W(\R)$ acts also on some space ${\bf C}$ of functions
on ${\open R}_N$.
Let us assume for definiteness that ${\bf C}$ is constituted by all
infinitely smooth functions of polynomial growth
\ie that $\forall f\in{\bf C}$ $\exists n(f)\in {\open Z}_+$ such that
$lim_{\x\to \infty} f(\x)/|\x|^{n(f)} =0$.
Every function $f\in{\bf C}$
can be also considered as element of $End\,{\bf C}$ \ie as operator on
${\bf C}$ acting by multiplication: $g\to fg$ $\forall g\in{\bf C}$.
The action of   $R_\vv\in W(\R)$ on such functions has the following form
\be
R_\vv (f(\x)) =f(R_\vv(\x)) \mbox{ when } f(\x)\in {\bf C},\\
R_\vv f(\x) =f(R_\vv(\x))R_\vv \mbox{ when } f(\x)\in End{\bf C}.
\ee

%by the rule $R_\vv |f(\x)\rangle=|f(\hat R_\vv(\x))\rangle$.

Dunkl differential-difference operators is defined as \cite{Dunkl}
\be\label{Dun}
D_i= \di +\sum_{\vv\in\R} \nu_\vv\frac {v_i} {(\x,\,\vv)} (1-R_\vv),
\ee
where coupling constants $\nu_\vv$ are such that
\be\label{nu}
\nu_\vv=\nu_{R_{\vec u}(\vv)} \ \forall \vec u,\, \vv\in\R.
\ee
These operators are well defined on ${\bf C}$ and
commute with each other, $[D_i,\, D_j]=0$.

Due to this property the deformed creation and annihilation operators
\cite{Poly, BHV}
\be\label{aa}
a_i^\alpha =\frac 1 {\sqrt{2}} (x_i + (-1)^\alpha D_i),\quad \alpha =0,1,
\ee
transform under action of reflections $R_\vv$ $\forall \vv\in\R$ as
${\open R}_N$ vectors
\be\label{comav}
R_\vv a_i^\alpha = \sum_{j=1}^N \left(\delta_{ij} - 2
\frac {v_i v_j}{(\vv,\,\vv)}\right)a_j^\alpha  R_\vv
\ee
and
satisfy the following commutation relations
\be\label{comaa}
[a^\alpha_i, a^\beta_j] = \varepsilon^{\alpha\beta}
\left(\delta_{ij}+
2\sum_{\vv\in\R} \nu_\vv \frac {v_i v_j}{(\vv,\,\vv)}R_\vv\right)
\ee
where $\varepsilon^{\alpha\beta}$
is the antisymmetric tensor, $\varepsilon^{01}=1$.

The operators $a_i^\alpha$ ($\alpha=0,1$, $i=1,...,N$)
together with the elements of the
group $W({\bf R})$ generate the associative algebra
with elements polynomial on
$a_i^\alpha$. We denote this algebra as $H_{W({\bf R})}(\nu)$ and
call it {\it the algebra of observables of Calogero model based on the
root system $\R$}. Here the notation
$\nu$ stands for a complete
set of $\nu_\vv$ with $\vv\in\R$.

The commutation relations (\ref{comaa}) and (\ref{comav}) allows one
to define the parity $\pi$:
\be\label{pi}
\pi (a_i^\alpha)=1\ \forall \alpha,\,i,
\qquad \pi(g)=0 \ \forall g\in W(\R)
\ee
and consider
$H_{W({\bf R})}(\nu)$ as a superalgebra.

Obviously the superalgebra $H_{W({\bf R})}(\nu)$ containes as a subalgebra
the group algebra ${\open C}[W(\R)]$ of Coxeter group $W(\R)$.

An important property of superalgebra $H_{W({\bf R})}(\nu)$ is that it
has ${\fr sl}_2$ algebra of inner differentiatings with the generators
\be\label{sl2}
T^{\alpha\beta}=\half \sum_{i=1}^N \left\{a_i^\alpha,\,a_i^\beta\right\}
\ee
which commute with ${\open C}[W(\R)]$, $[T^{\alpha\beta},\,R_\vv]=0$,
and act on $a_i^\alpha$ as on ${\fr sl}_2$-vectors:
\be\label{sl2vec}
\left[T^{\alpha\beta},\,a_i^\gamma\right]=
\varepsilon^{\alpha\gamma}a_i^\beta +\varepsilon^{\beta\gamma}a_i^\alpha.
\ee

The restriction of operator
$\hat H_{Cal}^\R(\nu)\defeq T^{01}$ on the subspace
${\bf C}_{W(\R)}\subset {\bf C}$ of $W(\R)$-invariant functions
is the second-order differential operator which is well-known
Hamiltonian of the rational Calogero model \cite{Cal} based on the
root system $\R$ \cite{OP}.
One of the relations (\ref{sl2}) namely $[\hat H_{Cal}^\R(\nu),\,a_i^\alpha]=
-(-1)^\alpha a_i^\alpha$ allows one to find the wavefunctions of
the equation $\hat H_{Cal}^\R(\nu)\psi =\epsilon \psi$
via usual Fock procedure with
the vacuum $ |0\rangle$ such that $a_i^0|0\rangle$=0 $\forall i$ \cite{BHV}.
After $W(\R)$-symmetrization these wavefunctions become the wavefunctions
of Calogero Hamiltonian $H_{Cal}^\R(\nu)|_{{\bf C}_{W(\R)}}$.

\section{Supertraces on $H_{W({\bf R})}(\nu)$}\label{trace}

{\sf\bf Definition}.
{\it The supertrace on
the superalgebra
%$H_{W(\R)}(\nu)$
${\cal A}$
is a linear
complex-valued function $str(\cdot)$ such that
$\forall f,g \in
{\cal A}$ with definite parity $\pi(f)$ and $\pi(g)$
\bee\label{scom}
str(fg)&=&(-1)^{\pi(f)\pi(g)}str(gf).
\eee
}

\vskip 2mm
Every supertrace $str(\cdot)$ on ${\cal A}$ generates
the invariant bilinear form on ${\cal A}$
\bee\label{bf}
B_{str}(f,g)=str(f\cdot g).
\eee

It is obvious that if such a bilinear form is degenerated then the
null-vectors of this form constitute both-sided
ideal ${\cal I}\subset {\cal A}$.

It was shown that the ideals of this sort are present
in the superalgebras $H_{W(A_1)}(\nu)$
(corresponding to the two-particle Calogero model) at $\nu =k+\half$ \cite{V}
and
in the superalgebras $H_{W(A_2)}(\nu)$
(corresponding to three-particle Calogero model)
at $\nu =k+\half$
and $\nu=k\pm\frac1 3$ \cite{K2} with every integer $k$ and that for all the
other values of $\nu$ all supertraces on these superalgebras
generate the nondegenerated bilinear forms (\ref{bf}).

The spectrum of $N$-particle rational Calogero Hamiltonian
(case ${\bf R}=A_{N-1}$) coincides with the spectrum
of system of $N$ noninteracting oscillators if the latter is shifted on
the constant $c(\nu) =\half N(N-1)\nu$.
It allows one to construct the similarity transformation between operators
$\hat H_{Cal}^{A_{N-1}}(\nu)|_{{\bf C}_{W(\R)}} - c(\nu)$
with different $\nu$ \cite{GP}.
Nevertheless the previous consideration shows that
the corresponding algebras $H_{W(A_{N-1})}(\nu)$ can be nonisomorphic
at different values of $\nu$.

It is easy to describe all supertraces on
${\open C}[W(\R)]$. Every supertrace on ${\open C}[W(\R)]$
is completely defined by its values on $W(\R)\subset {\open C}[W(\R)]$
and the function $str$ is a central function on $W(\R)$ \ie
the function on the conjugacy classes.

To formulate the theorem establishing the connection between
the supertraces on $H_{W(A_N)}(\nu)$ and the supertraces on
${\open C}[W(\R)]$
let us introduce the grading on the vector space of ${\open C}[W(\R)]$.
The grading $E$ of elements $g\in W(\R)\subset
{\open C}[W(\R)]$ is defined as follows.
Let ${\cal H}_N^\alpha$ be the linear space with basis elements
$a_1^\alpha,\ a_2^\alpha,\,...\,,\,a_N^\alpha$.
Consider the subspaces
${\cal E}^\alpha (g) \subset {\cal H}_N^\alpha$ as
\be\label{Halpha}
{\cal E}^\alpha (g) =\{h\in {\cal
H}_N^\alpha:\quad gh=-hg \}\,,
\ee
and put
\be\label{grad}
E(g)=dim\,{\cal E}^\alpha(g).
\ee
To avoid misunderstanding it should be noticed that
${\open C}[W(\R)]$ is not in general a graded algebra.

The following theorem was proved in \cite{KV}
\footnote{this theorem was proved for the case $\R=A_N$ only
but the proof does not depend on the particular properties
of the symmetric group $S_N=W(A_{N-1})$.}:

{\bf Theorem 1.} {\it Let ${\cal P}(g)$ be the projector
${\open C}[W(\R)] \rightarrow {\open C}[W(\R)]$
defined as ${\cal P}(\sum_i \alpha_i
g_i)=\sum_{i:\, g_i \neq {\bf 1}} \alpha_i g_i$  ($g_i\in W(\R)$,
$\alpha_i\in {\open C})$. Let the grading $E$ defined in (\ref{grad})
and the subspaces ${\cal E}^\alpha(g)$ defined in (\ref{Halpha})
satisfy the equations
\bee\label{main}
E({\cal P}([h_0,\,h_1])g)=E(g)-1 \ \ \forall g\in {\open C}[W(\R)], \ \
\forall h_\alpha\in {\cal E}^\alpha(g).
\eee
Then every supertrace on the algebra ${\open C}[W(\R)]$
satisfying the equations
\bee\label{GLC}
str([h_0,\,h_1]g)=0\qquad \forall g\in W(\R) \mbox{
with }E(g)\neq 0 \mbox{ and }\forall h_\alpha \in {\cal E}^\alpha(g),
\eee
can be extended to the supertrace on $H_{W(\R)}(\nu)$ in a unique way.}

%%%%%%%%%%%%%%%%%%%%%%%%%%%%%%%%%%%%%%%%%%%%%%%%%%%%%%%%%%%%%%%%%%%%%%%%%%

\section{Superalgebras $H_{W(I_2(n))}$}\label{I2}

It is convenient to use ${\open C}$ instead of ${\open R}_2$
to describe $H_{W(I_2(n))}$.
The root system $I_2(n)$ contains $2n$ vectors $v_k=exp(\pi i k/n)$,
$k=0,1,...,2n-1$. The corresponding Coxeter group $W(I_2(n))$
has $2n$ elements, $n$ reflections $R_k$ acting on $z,\ z^* \in{\open C}$
as follows
\bee\label{Rk}
R_k z&=& -z^* v_k^2 R_k, \nn
R_k z^* &=& -z {v_k^*}^2 R_k,\qquad k\in {\open Z}_n
\eee
and $n$ elements of the form $S_k=R_k R_0$. $S_0$ is the unity in $W(I_2(n))$.
These elements satisfy the following relations
\be\label{RR}
R_k R_l = S_{k-l},\qquad
S_k S_l = S_{k+l},\qquad
R_k S_l = R_{k-l},\qquad
S_k R_l = R_{k+l}.
\ee

Obviously the reflections $R_{2k}$ lie in one conjugacy class
and $R_{2k+1}$ in another if $n$ is even. If $n$ is odd then all reflections
$R_k$ lie in one conjugacy class.

It is convenient to consider the following basis in ${\open C}[W(I_2(n))]$
\be\label{LQbas}
L_p=\sum_{k=0}^{n-1} \lambda^{kp} R_k,\qquad
Q_p=\sum_{k=0}^{n-1} \lambda^{-kp} S_k,\\
\mbox{where }\lambda=exp(2\pi i/n),\qquad p\in{\open Z}_n.
\ee

Differential-difference operators have the following form
\bee\label{Dz}
D_{z}=2\frac \partial {\partial z^*} +
   \nu_0 \sum_{k=0}^{n-1} \frac {v_{2k}} {zv_{2k}^* + z^*v_{2k}}
   \left(1-R_{2k}\right)
+   \nu_1 \sum_{k=0}^{n-1} \frac {v_{2k+1}} {zv_{2k+1}^* + z^*v_{2k+1}}
   \left(1-R_{2k+1}\right),         \nn
D_{z^*}=2\frac \partial {\partial z} +
   \nu_0 \sum_{k=0}^{n-1} \frac {v_{2k}^*} {zv_{2k}^* + z^*v_{2k}}
   \left(1-R_{2k}\right)
+   \nu_1 \sum_{k=0}^{n-1} \frac {v_{2k+1}^*} {zv_{2k+1}^* + z^*v_{2k+1}}
   \left(1-R_{2k+1}\right).
\eee

This form unifies both cases of even and odd $n$. If $n$ is odd then $D_z$
and $D_{z^*}$ depends on the only coupling constant $\nu_0+\nu_1$.

The basis in the space ${\cal H}_2^\alpha$ is
\bee\label{ab}
&{}& a^\alpha \defeq \half(z   + (-1)^\alpha D_z    ) \nn
&{}& b^\alpha \defeq \half(z^* + (-1)^\alpha D_{z^*}).
\eee

Now we can write down the relation between elements $a^\alpha$, $b^\alpha$,
$R_k$ and $S_k$ generating the associative superalgebra $H_{W(I_2(n))}$
of polynomials of $a_i^\alpha$, $b_i^\alpha$ with coefficients in
${\open C}[W(I_2(n))]$:
\be
% G-H
\label{RH}
R_k a^\alpha=  - {v_k  }^2 b^\alpha R_k ,&\qquad&
R_k b^\alpha = - {v_k^*}^2 a^\alpha R_k , \\
\label{SH}
S_k a^\alpha =    {v_k^*}^2 a^\alpha S_k ,&\qquad&
S_k b^\alpha =   v_k^2     b^\alpha S_k , \\
\label{LH}
L_p a^\alpha = -   b^\alpha L_{p+1} ,&\qquad&
L_p b^\alpha = -   a^\alpha L_{p-1} ,\\
\label{QH}
Q_p a^\alpha = a^\alpha Q_{p+1} ,&\qquad&
Q_p b^\alpha = b^\alpha Q_{p-1} ,\\
% G-G
\label{GG}
L_k L_l = n \delta_{k+l} Q_l,&\qquad& L_k Q_l = n \delta_{k-l} L_l,\nn
Q_k L_l = n \delta_{k+l} L_l,&\qquad& Q_k Q_l = n \delta_{k-l} Q_l,
\ \ \mbox{where }\delta_k\defeq \delta_{k0},
\ee
% H-H
\bee\label{HH}
\left[a^\alpha,\,b^\beta\right]&=&\varepsilon^{\alpha\beta}
\left( 1+\frac {\nu_0+\nu_1}2 L_0 -\frac {\nu_0-\nu_1}2 L_{n/2}\right), \nn
\left[a^\alpha,\,a^\beta\right]&=&\varepsilon^{\alpha\beta}
\left( \frac {\nu_0+\nu_1}2 L_1 -\frac {\nu_0-\nu_1}2 L_{n/2+1} \right),\nn
\left[b^\alpha,\,b^\beta\right]&=&\varepsilon^{\alpha\beta}
\left( \frac {\nu_0+\nu_1}2 L_{-1} -\frac {\nu_0-\nu_1}2 L_{n/2-1} \right).
\eee
The terms containing $\nu_0-\nu_1$ in (\ref{HH}) are absent when $n$ is odd
and halfinteger indices are senseless, so let us assume that
$\nu_0$=$\nu_1$ when $n$ is odd.
% and do not stipulate every time if $n$ is even or not.

\section{The number of supertraces on $H_{W(I_2(n))}$}
In this section the following theorem is proved:

{\bf Theorem 2}. {\it The superalgebra $H_{W(I_2(n))}$
has $\left[\frac {n+1} 2\right]$ supertraces.}

\vskip 1mm
\noindent
It is easy to find the grading $E$:
\be\label{E}
&{}& E(R_k)=1,\qquad R_k (v_k^* a^\alpha + v_k b^\alpha) =
- (v_k^* a^\alpha + v_k b^\alpha) R_k,          \nn
&{}& E(S_k)=0,\qquad \forall k\neq \frac n 2  \,,  \\
&{}& E(S_{n/2})=2,\qquad
S_{n/2}a^\alpha =-a^\alpha  S_{n/2},\ S_{n/2}b^\alpha =-b^\alpha  S_{n/2},
\mbox{ if $n$ is even and $S_{n/2}$ exists,}   \nonumber
\ee
and check that  $H_{W(I_2(n))}$
satisfies the conditions of Theorem 1, \ie that
\be\label{E1}
E\left({\cal P}\left(\left[(v_k^* a^0 + v_k b^0),\,
(v_k^* a^1 + v_k b^1)\right]\right) R_k\right)=0
\ee
and that for even $n$
\be\label{E2}
&{}& E\left({\cal P}\left(\left[a^0,\,
a^1 \right]\right) S_{n/2}\right)=
E\left({\cal P}\left(\left[a^0,\,
b^1 \right]\right) S_{n/2}\right)=\nn
&{}& E\left({\cal P}\left(\left[b^0,\,
a^1 \right]\right) S_{n/2}\right)=
E\left({\cal P}\left(\left[b^0,\,
b^1 \right]\right) S_{n/2}\right)=1.
\ee

To compute the number of supertraces on the superalgebra
we have to find the number of the solutions of the equations (\ref{GLC})
which have the following form for the algebra under consideration
\be\label{GLC1}
str\left(\left[(v_k^* a^0 + v_k b^0),\,
(v_k^* a^1 + v_k b^1)\right] R_k\right)=0
\ee
and
\be\label{GLC2}
&{}& str\left(\left[a^0,\,a^1 \right] S_{n/2}\right)=
str\left(\left[a^0,\,b^1 \right] S_{n/2}\right)=\nn
&{}& str\left(\left[b^0,\,a^1 \right] S_{n/2}\right)=
str\left(\left[b^0,\,b^1 \right] S_{n/2}\right)=0
\qquad\mbox{for even $n$}.
\ee
The equations  (\ref{GLC1}) lead to
\be\label{strT}
str(R_k)=&-&\frac {\nu_0+\nu_1} 2
str \left( Q_0+\half (Q_{1}+Q_{-1})\right)\nn
   &+&(-1)^k\,\,
\frac {\nu_0-\nu_1} 2 str \left( Q_{n/2}+\half (Q_{n/2+1}+Q_{n/2-1})\right),
\ee
and when $n$ is even the relations (\ref{GLC2}) take place and lead  to
\be\label{strS1}
&{}& str \left( (\nu_0 + \nu_1 ) L_{\pm 1}- (-1)^{n/2}
(\nu_0 - \nu_1) L_{n/2 \pm 1} \right)=0   \\
\label{strS2}
&{}& str(S_{n/2})=-\half
str \left( (\nu_0 + \nu_1 ) L_{0}- (-1)^{n/2}
(\nu_0 - \nu_1) L_{n/2} \right).
\ee
It is easy to see that equations (\ref{strS1}) are consequences
of (\ref{strT}).
The equations (\ref{strT}) express $str(R_k)$ $\forall k$ via $str(S_j)$
and (\ref{strS2}) expresses $str(S_{n/2})$ via $str(S_j)$ with $j\neq n/2$
when $n$ is even.
Hence due to theorem 1 every supertrace on $H_{W(I_1(n))}(\nu)$ is
determined completely by its values on $S_j$ with $j\neq n/2$.
Since $S_j$ and $S_{k}$ belong to one conjugacy class if and only if
$j=-k$ the number of independent supertrace is equal to $(n+1)/2$
when $n$ is odd and to $n/2$ if $n$ is even.
It finishes the proof of Theorem 2.

For $n=2,3,4$ the obtained result are in agreement with the
results obtained in \cite{KV} and \cite{K}. The case $n=6$
gives that the superalgebra of observables
of the rational Calogero model based on the root system of
$G_2$ type has 3 independent supertraces.

\end{document}